\documentclass{article}

\usepackage{arxiv}

\usepackage[utf8]{inputenc} 
\usepackage[T1]{fontenc}    
\usepackage{hyperref}       
\usepackage{url}            
\usepackage{booktabs}       
\usepackage{amsfonts}       
\usepackage{nicefrac}       
\usepackage{microtype}      
\usepackage{lipsum}
\usepackage{graphicx}
\graphicspath{ {./images/} }
\usepackage{lipsum}
\usepackage[round, authoryear]{natbib}
\usepackage{slashbox}
\usepackage{color}
\usepackage{tikz}

\title{Neural networks to predict survival from RNA-seq data in oncology}
\author{
  Mathilde Sautreuil \\
Laboratoire MICS, CentraleSupélec, Université Paris-Saclay \\
9 rue Joliot Curie, 91190 Gif-sur-Yvette, France \\ 
\texttt{prenom.nom@centralesupelec.fr} \\
\AND
Sarah Lemler \\
Laboratoire MICS, CentraleSupélec, Université Paris-Saclay \\
9 rue Joliot Curie, 91190 Gif-sur-Yvette, France \\ 
\texttt{prenom.nom@centralesupelec.fr} \\
\And
Paul-Henry Cournède \\
Laboratoire MICS, CentraleSupélec, Université Paris-Saclay \\
9 rue Joliot Curie, 91190 Gif-sur-Yvette, France \\ 
\texttt{prenom.nom@centralesupelec.fr} \\
}

\begin{document}
\maketitle

\begin{abstract}
Survival analysis consists of studying the elapsed time until an event of interest, such as the death or recovery of a patient in medical studies. This work explores the potential of neural networks in survival analysis from clinical and RNA-seq data. If the neural network approach is not recent in survival analysis, methods were classically considered for low-dimensional input data. But with the emergence of high-throughput sequencing data, the number of covariates of interest has become very large, with new statistical issues to consider. We present and test a few recent neural network approaches for survival analysis adapted to high-dimensional inputs.
\end{abstract}


\section{Introduction}
\label{sec:Intro}
Survival analysis consists of studying the elapsed time until an event of interest, such as the death or recovery of a patient in medical studies.
This paper aims to compare methods to predict a patient's survival from clinical and gene expression data. 

The Cox model \citep{coxregression1972} is the reference model in the field of survival analysis. It relates the survival duration of an individual to the set of explanatory covariates. It also enables to take into account censored data that are often present in clinical studies. With high-throughput sequencing techniques, transcriptomics data are more and more often used as covariates in survival analysis. Adding these covariates raise issues of high-dimensional statistics, when we have more covariates than individuals in the sample. 
Methods based on regularization or screening \citep{tibshiranilasso1997,fanhigh-dimensional2010} have been developed and used to solve this issue. 


The Cox model relies on the proportional hazard hypothesis, and in its classical version, does not account for nonlinear effects or interactions, which proves limited in some real situations. Therefore, in this paper, we focus on another type of methods: neural networks. Deep learning methods are more and more popular, notably due to their flexibility and their ability to handle interactions and nonlinear effects, including in the biomedical field~\citep{rajkomarscalable2018,kwongclinical2017,suodeep2018}. 
The use of neural networks for survival analysis is not recent, since it dates back to the 90's \citep{faraggineural1995,biganzolifeed1998} but it began being widely used only recently. We can differentiate two strategies. The first one relies on the use of a neural network based on the Cox partial log-likelihood as those developed by \cite{faraggineural1995,chingcox-nnet:2018,katzmandeepsurv:2018,kvamme2019}. The second strategy consists of using a neural network based on a discrete-time survival model, as introduced by \cite{biganzolifeed1998}.
\cite{biganzolifeed1998} have studied this neural network only in low-dimension. In this paper, our objective is to study and adapt this model to the high-dimensional cases, and compare its performances to two other methods:
 the two-step procedure with the classical estimation of the parameters of the Cox model with a Lasso penalty to estimate the regression parameter and a kernel estimator of the baseline function (as in \cite{guilloux2016adaptive}) and the Cox-nnet neural network \citep{chingcox-nnet:2018} based on the partial likelihood of the Cox model.
 Section 2 recalls the different notations used in survival analysis and presents the different models. Then, we introduce the simulation plan created to compare the models. Finally, we underline the results to conclude with the potential of neural networks in survival analysis.


\section{Models}
\label{sec:Models}
First, we introduce the following notations:
\begin{itemize}
    \item $Y_i$ the survival time
    \item $C_i$ the censorship time 
    \item $T_i = \min(Y_i, C_i)$ the observed time
    \item $\delta_i$ the censorship indicator (which will be equal to 1 if the interest event occurs and else to 0).
\end{itemize}

\subsection{The Cox model}
\label{sec:Cox}
The Cox model~\citep{coxregression1972} predicts the survival probability of an individual from explanatory covariates $X_{i.} = (X_{i1} ,\dots, X_{ip})^{T} \in \mathbb{R}^p$. The hazard function $\lambda$ is given by:
\begin{equation}
    \lambda(t|X_{i.}) = \alpha_0(t)\exp(\beta^TX_{i.}),
    \label{eq:Cox}
\end{equation}
where $\alpha_0(t)$ corresponds to the baseline hazard and $\beta =(\beta_1,\dots,\beta_p)^T\in\mathbb{R}^p$ is the vector of regression coefficients.
A benefit of this model is that only $\alpha_0(t)$ depends on time while the second term of the right hand side of (\ref{eq:Cox}) depends only on the covariates (proportional hazard model). The Cox model structure can be helpful when we are interested in the prognostic factors because $\beta$ can be estimated without knowing the function $\alpha_0$. It is possible thanks to the Cox partial log-likelihood, which is the part of the total log-likelihood that does not depend on $\alpha_0(t)$, and is defined by:
$$\mathcal{L}\left(\beta\right)=\sum_{i=1}^n \left(\beta^TX_{i.}  \right) - \sum_{i=1}^{n} \delta_i \log\left(\sum_{l\in R_i} \exp\left( \beta^TX_{l.} \right)\right),$$
with $R_i$ the individuals at risk at the observation time $T_i$ of individual $i$, and $\delta_i$ the censorship indicator of individual $i$. 
The Lasso procedure was proposed by \cite{tibshiranilasso1997} for the estimation of $\beta$ in the high-dimensional setting. The non-relevant variables are set to zero thanks to the $L_1-$ penalty added to the Cox partial likelihood: $\mathcal{L}(\beta) + \lambda||\beta||_1.$
However, to predict the survival function $S$, we need to fully estimate the hazard risk $\lambda(s|X_{i.})$ since:
$$S(t) = \mathbb{P}(T_i > t | X_{i.}) = \exp\left(- \int_0^t \alpha_0(s) \exp(\beta^TX_{i.} )ds \right).$$
We follow the two-step procedure of \cite{guilloux2016adaptive}: first, we estimate $\beta$ from the penalized Cox partial likelihood, and then we estimate $\alpha_0(t)$ from the kernel estimator introduced by \cite{ramlau-hansensmoothing1983}, in which we have plugged the Lasso estimate of $\beta$.

\subsection{Neural networks}\label{sec:NN}
The studied neural networks in this paper are fully-connected multi-layer perceptrons. Several layers constitute this network with at least one input layer, one output layer, and one or several hidden layers. 
\subsubsection{Cox-nnet}
\label{subsec:chap3_coxnn}
\cite{faraggineural1995} developed a neural network based on the proportional hazards model. The idea of \cite{faraggineural1995} was to replace the linear prediction of the Cox regression with the neural network's hidden layer's output.
\cite{faraggineural1995} only applied their neural network to survival analysis from clinical data, in low dimension. More recently, some authors revisited this method \cite{chingcox-nnet:2018,katzmandeepsurv:2018,kvamme2019}. However, only Cox-nnet \citep{chingcox-nnet:2018} was applied in a high-dimensional setting. We will thus use this model as benchmark in our study.

The principle of Cox-nnet is that its output layer corresponds to a Cox regression: the output of the hidden layer replaces the linear function of the covariates in the exponential of the Cox model equation.

To estimate the neural network weights, \cite{chingcox-nnet:2018} uses the Cox partial log-likelihood as the neural network loss: 
\begin{equation}
    \mathcal{L}(\beta, W, b) = \sum_{i=1}^n \theta_i - \sum_{i=1}^{n} \delta_i \log\left(\sum_{l\in R_i} \exp\left(\theta_l\right)\right)
\end{equation}
with $\delta_i$ the censoring indicator and $\theta_i =\beta^TG(W^TX_{i.}+b)$, where $G$ is the activation function of the hidden layer, $W = (w_{dh})_{1\le d \le p, 1 \le h \le H}$ with $H$ the number of neurons in the hidden layer, and $\beta = (\beta_1, \ldots, \beta_H)^T$ the weights and $b$ the biases of the neural network to be estimated. In this network, the activation function $tanh$ is used. To the partial log-likelihood, \cite{chingcox-nnet:2018} adds a ridge penalty in $L_2-$norm of the parameters. Thus, the final cost function for this neural network is: 
\begin{equation}
Loss(\beta, W, b) = \mathcal{L}(\beta,W,b) + \lambda(\|\beta\|_2 + \|W\|_2 + \|b\|_2). \label{eq:Cout_coxnn}    
\end{equation}

We maximize this loss function to deduce estimators of $\beta$, $W$ and $b$. The principle in this neural network is that the activation function for the output layer is a Cox regression, so that we have: 
\begin{equation}
\label{eq:h}
    \hat{h}_i = \exp \left( \underbrace{\sum_{h=1}^{H} \hat\beta_h G \left(\hat b_h + \hat W^TX_{i.}\right)}_{\hat\theta_i = \hat\beta^TG(\hat W^TX_{i.}+\hat b)} \right). 
\end{equation}
The output of the neural network $\hat{h}_i$ corresponds to the part of the Cox regression that does not depend on time. \cite{chingcox-nnet:2018} only used $\hat{h}_i$, but in our study, we are interested in the complete survival function, and thus we need to estimate the complete hazard function $\hat{h}(x_i,t)$. For that purpose, we estimate the baseline risk $\alpha_0(t)$, with the kernel estimator introduced by \cite{ramlau-hansensmoothing1983}.
As for the Cox model, we estimate $\alpha_0(t)$ with the two-steps procedure of \cite{guilloux2016adaptive} and this estimator is defined by: 
\begin{equation}
\widehat{\alpha}_{m}(t) = \frac{1}{n m} \sum_{i=1}^{n} K \left(\frac{t-u}{m} \right) \frac{\delta_i}{\sum_{l \in R_i} \widehat{h}_l},
    \label{eq:RH}
\end{equation}
with $\hat{h}_l$ the estimator defined by (\ref{eq:h}), $K: \mathbb{R} \rightarrow \mathbb{R}$ a kernel (a positive function with integral equal to 1), $m$ the bandwidth, which is a strictly positive real parameter. $m$ can be obtained by cross-validation or by the Goldenshluger \& Lepski method~\cite{goldenshluger_bandwidth_2011} for instance, and we choose the latter.
We can finally derive an estimator of the survival function for individual $i$: 
\begin{equation}
    \widehat{S}(t|X_{i.}) = \exp \left( - \int_{0}^t \widehat{\alpha}_{m}(s) \hat h_i ds\right).
\end{equation}

\subsubsection{Discrete time neural network}
\label{subsec:chap3_NNsurv1}
\begin{figure}
\includegraphics[width=\textwidth]{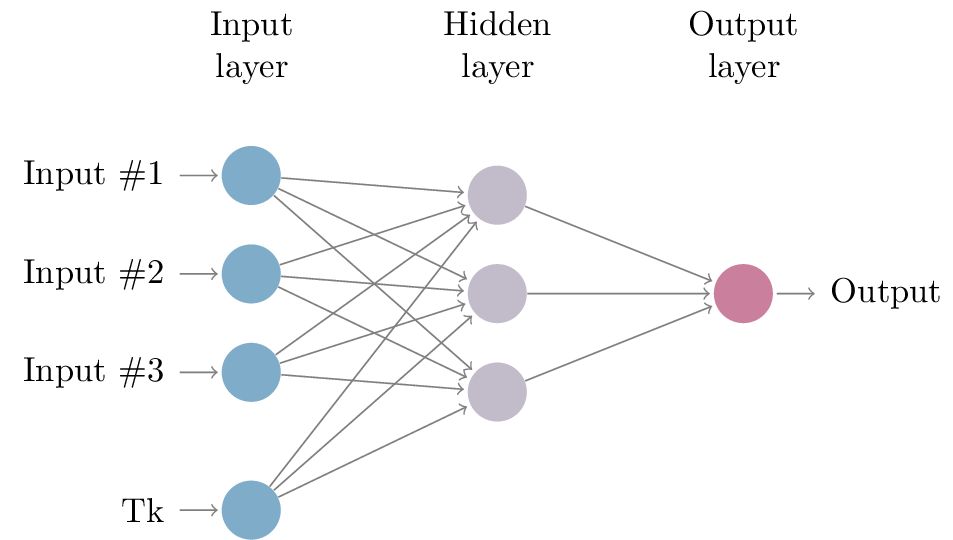}
    \caption{Structure of the neural network based on the discrete-time model of \cite{biganzolifeed1998}}
    \label{fig:NNsurv}
\end{figure}

\cite{biganzolifeed1998} have proposed a neural network based on a discrete-time model. They introduced $L$ time intervals $A_l = ]t_{l-1},t_l]$, and build a model predicting in which interval, the failure event occurs. We write the discrete hazard as the conditional probability of survival: 
\begin{equation}
    h_{il} = P(Y_i \in A_l|Y_i> t_{l-1}), \label{eq:hazDis}
\end{equation}
with $Y_i$ the survival time of individual $i$. \cite{biganzolifeed1998} duplicates the individuals as input of the neural network. The duplication of individuals gives it a more original structure than that of a classical multi-layer perceptron. The \cite{biganzolifeed1998}'s neural network takes as input the set of variables of the individual and an additional variable corresponding to the mid-point of each interval. Due to the addition of this variable, the $p$ variables of each individual are repeated for each time interval. The output is thus the estimated hazard ${h}_{il} = h_l(X_i, a_l)$ for the individual $i$ at time $a_l.$ We schematize the structure of this neural network on \textsc{Figure}~\ref{fig:NNsurv}.
\cite{biganzolifeed1998} initially used a 3-layers neural network with a logistic function as the activation function for both the hidden and output layers. 
The output of the neural network with $H$ neurons in the hidden layer and $p+1$ input variables is given by: 
$${h}_{il} = {h}(x_i, t_l) = f_2\left(a + \beta^T f_1\left(b + W^T X_{i.} \right)\right),$$
where  $W = (w_{dh})_{1\le d \le p+1, 1 \le h \le H}$, and $\beta = (\beta_{1}, \ldots, \beta_{H})^T$ are the weights of the neural network, $a$ and $b$ are the biases of the neural network to be estimated, and $f_1$ and $f_2$ the sigmoid activation functions. The target of this neural network is the death indicator $d_{il}$, which will indicate if the individual $i$ dies in the interval $A_l$. We introduce $l_i \le L$ the number of intervals in which the individual $i$ is observed, $d_{i0},\ldots,d_{i(l_i-1)} = 0$ whatever the status of the individual $i$ and $d_{il_i}$ is equal to $0$ if the individual $i$ is censored and $1$ otherwise.
The cost function used by \cite{biganzolifeed1998} is the cross-entropy function and the weights of the neural network can be estimated by minimizing it: 
\begin{eqnarray}
    \mathcal{L}(\beta,W,a,b) = -\sum_{i=1}^{n} \sum_{l=1}^{l_i} d_{il} \log({h}_{il}) + (1-d_{il}) \log(1 - {h}_{il}). \label{eq:crossEntropy}
\end{eqnarray}

The duplication of the individuals for each time interval increases the sample size in the neural network, it is an advantage in a high-dimensional framework. Moreover, \cite{biganzolifeed1998} added a ridge penalty to their cross-entropy function (\ref{eq:crossEntropy}): 
\begin{equation}
\label{eq:loss}
    Loss(\beta,W,a,b) = \mathcal{L}(\beta,W,a,b)+ \lambda (\|\beta\|_2+\|W\|_2+\|a\|_2+\|b\|_2),
\end{equation}
In \cite{biganzolifeed1998}, $\lambda$ was chosen by deriving an Information Criteria. We choose instead to use cross-validation since it improved model the predictive capacity.

    
After estimating the parameters of the neural network by minimizing the loss function (\ref{eq:loss}), the output obtained is the estimate of the discrete risk $\widehat{h}_{il}$ for each individual $i$ and the survival function of individual $i$ is estimated using: 
\begin{equation}
     \widehat{S}(T_{l_i}) = \prod_{i=1}^{l_i} (1 - \widehat{h}_{il}). \label{eq:hazTosurv}
\end{equation}
This model was only applied for low-dimensional inputs, and this paper investigates its performance and capacity to adapt to high-dimensional settings. We denote this network NNsurv. We noticed an improvement of the performance when using a ReLU activation function for the hidden layers and thus used it instead of the original sigmoid functions. Moreover, original neural network only has one hidden layer. We propose to add one supplementary hidden layer to study if a deeper structure could improve the neural network prediction capacity. We call the deeper version NNsurv-deep. Its structure is similar to the one schematized in Figure~\ref{fig:NNsurv}, but with two hidden layers instead of one. The input layer does not change, and the individuals are always duplicated at the input of the neural network. The output layer also has a single neuron corresponding to the discrete hazard estimate. These neural networks are implemented in a package available on \url{https://github.com/mathildesautreuil/NNsurv}.

We will compare the performances of these four models (Cox-Lasso, Cox-nnet, NNsurv, NNsurv-deep) on simulated data and then to a real dataset.


\section{Simulations}
\label{sec:Simu}
We create a simulation design to compare different neural network approaches to predict survival time in high-dimension. We divide the simulation plan into two parts. The first part concerns a simulation study based on \citep{bendergenerating2005} which proposes to generate the survival data from a Cox model. Data simulated with this model naturally favors the two methods based on the Cox model. We also consider a model with a more complex behavior: the Accelerated Hazards (AH) model~\citep{chen_analysis_2000}. In the AH model, variables will accelerate or decelerate the hazard risk. The survival curves of the AH model can therefore cross each other. Other choices of models were also possible, and in the Appendix \ref{appendix}, we also present the results for the Accelerated Failure Time (AFT) model~\citep{kalbfleischstatistical2002} which does not satisfy the proportional risk assumption either, but does not allow the intersection of survival curves of different patients. 

In all cases, the models' baseline risk function is assumed known and follows a particular probability distribution. We use the Weibull distribution for the Cox model and the log-normal distribution for the AH model. Several simulations are considered, by varying the sample size, the total number of explanatory variables, and the number of relevant explanatory variables considered in the model. We use the package that we have developed called \verb?survMS? and available on CRAN or \url{https://github.com/mathildesautreuil/survMS}.

\subsection{Generation of survival times}
Considering the survival models (Cox, AFT, and AH models), the survival function $S(t|X)$ can be written as: 
\begin{equation}
    S(t|X) = \exp(-H_0(\psi_1(X)t)\psi_2(X)\;\mbox{with}
\end{equation}
$H_0$ is the cumulative hazard and
 $$
(\psi_1(X), \psi_2(X)) = \left\{
    \begin{array}{ll}
        (1, \exp(\beta^TX)) & \mbox{for the Cox model } \\
        (\exp(\beta^TX), \exp(-\beta^TX)) & \mbox{for the AH model} \\
        (\exp(\beta^TX), 1) & \mbox{for the AFT model. } 
    \end{array}
\right.
$$

 The distribution function is deduced from the survival function from the following formula: 
 \begin{equation}
     F(t|X) = 1 - S(t|X). \label{eq:F}
 \end{equation}
 For data generation, if $Y$ is a random variable that follows a probability distribution $F,$ then $U = F(Y)$ follows a uniform distribution on the interval $[0,1],$ and $(1-U)$ also follows a uniform distribution $\mathcal{U}[0,1].$ From Equation (\ref{eq:F}), we finally obtain that: 
\begin{equation}
  1 - U = \exp(-H_0(\psi_1(X)t)\psi_2(X)).
\end{equation}

 If $\alpha_0(t)$ is positive for all $t,$ then $H_0(t)$ can be inverted, and we can express the survival time of each of the models considered (Cox, AFT and AH) from $H_0^{-1}(u).$ We write in a general form the expression of the random survival times for each of the survival models: 
 \begin{equation}
      T = \frac{1}{\psi_1(X)} H^{-1}_0 \left( \frac{\log(1-U)}{\psi_2(X)} \right).
  \end{equation}
Two distributions are used for the cumulative hazard function $H_0(t)$ to generate the survival data. If the survival times are distributed according to a Weibull distribution $\mathcal{W}(a, \lambda),$ the baseline hazard is of the form
: 
\begin{equation}
    \alpha_0(t) = a\lambda t^{a-1}, \lambda > 0, a > 0.
\end{equation}
 The inverse of the cumulative risk function is expressed as follows:
 \begin{equation}
     H_0^{-1}(u) = \left( \frac{u}{\lambda} \right)^{1/a}.
 \end{equation}
 
For survival times following a log-normal distribution $\mathcal{LN}(\mu, \sigma)$ with mean $\mu$ and standard deviation $\sigma$, the basic risk function is therefore written: 
 \begin{equation}
    \alpha_0(t) = \frac{\frac{1}{\sigma\sqrt{2\pi t}} \exp\left[-\frac{(\log t - \mu)^2 }{2 \sigma^2}\right]}{1 - \Phi\left[\frac{\log t - \mu}{\sigma}\right]},
 \end{equation}
 with $\Phi(t)$ the distribution function of a standard Normal distribution. 
 The inverse of the cumulative hazard function is expressed by: 
 \begin{equation}
     H_0^{-1}(u) = \exp(\sigma\Phi^{-1}(1-\exp(-u))+\mu), \label{eq:InvHazLN}
 \end{equation}
 with $\Phi^{-1}(t)$ the inverse of the distribution function of a centered and reduced normal distribution.
 

\subsection{Simulation with the Cox - Weibull model}


\paragraph{\textbf{Survival times and baseline function: }}
Generating survival times from a variety of parametric distributions were described by \cite{bendergenerating2005}.
In the case of a Cox model with a baseline function distributed from a Weibull distribution, the inverse cumulative hazard function is $H_0^{-1}(t) = (\frac{t}{\lambda})^{\frac{1}{a}}$ and the survival time T of the Cox model is expressed as:
\begin{equation}
    T = \left( -\frac{1}{\lambda} \log(1-U) \exp(-X_{i.} \beta) \right)^{\frac{1}{a}}, \label{eq:Times}
\end{equation}
where $U$ is a random variable with $U \sim \mathcal{U}[0,1].$




\paragraph{\textbf{Choice of parameters of the Weibull distribution:}}
 \label{sec:Param}
We chose the Weibull distribution parameters so that our design of simulation is close to real datasets. The mean and the standard deviation of Breast cancer real dataset (available on \url{www.ncbi.nlm.nih.gov/geo/query/acc.cgi?acc=GSE6532}) is around $2325$ days and $1304$ days respectively. As the survival times follow a Weibull distribution, the mean and the variance of $T$ write as:
$$\mathbb{E}(T) = \frac{1}{\sqrt[a]{\lambda}} \Gamma\left(\frac{1}{a}+1\right) \,\; \mbox{and} \,\; \mathbb{V}(T) = \frac{1}{\sqrt[a]{\lambda^2}} \left[ \Gamma\left(\frac{2}{a}+1\right) - \Gamma^2\left(\frac{1}{a}+1\right)\right],$$
where $\Gamma$ is the Gamma function. We set $a = 2$ and $\lambda = 1.3e^{-7}$ to have a mean and variance of our simulated datasets close to those of the Breast cancer real dataset.

 \subsection{Simulation with the AH - Log-Normal model}
\paragraph{\textbf{Survival times and baseline function: }}
 Building on the work of \cite{bendergenerating2005}, we also simulate the survival data from the AH model. We perform this simulation to generate data whose survival curves will intersect. 
For this simulation, we consider that the survival times follow a log-normal distribution $\mathcal{LN}(\mu,\sigma).$ In this case, the inverse of the cumulative hazard function is expressed as (\ref{eq:InvHazLN}), and we have: 
\begin{equation}
    T = \frac{1}{\exp(\beta^T X_{i.})} \sigma \Phi^{-1}\left(\frac{\log(1-U)}{\exp(-\beta^T X_{i.})} + \mu\right) \label{eq:Tah}
\end{equation}
 with $\Phi^{-1}(t)$ the inverse of the distribution function of a centered and reduced normal distribution.
 
\paragraph{\textbf{Choice of parameters of the Log-Normal distribution:}}
As in the previous simulation, we ensure that the distribution of the simulated data is close to that of the real ones and we use the formulas: 
\begin{equation}
\mu = \ln(\mathrm{E}(T))- \frac{1}{2} \sigma^2 \label{eq:mu} \,\;\, \mbox{and} \,\;\,
\sigma^2 = \ln\left(1+\frac{\mathrm{Var}(T)}{(\mathrm{E}(T))^2}\right).
\end{equation}
Since, the expectation and the standard deviation are respectively $2325$ and $1304$, the values of $\mu$ and $\sigma$ used for the simulation of the survival data should be $\mu = 7.73$ and $\sigma = 0.1760.$ However, to have survival curves crossing rapidly, we take a higher value of $\sigma$: $\sigma=0.7$.
\subsection{Metrics}
\label{subsec:Metric}
To assess the performance of the survival models, we use two classical metrics, the Concordance Index (CI) and the Integrated Brier Score (IBS).
\subsubsection{Concordance Index}
The index measures whether the prediction of the model under study matches the rank of the survival data. If the event time of an individual $i$ is shorter than that of an individual $j$, a good model will predict a higher probability of survival for individual $j$. This metric takes into account censored data, and it takes a value between $0$ and $1$. If the C-index is equal to $0.5$, the model is equivalent to random guessing.
The time-dependent C-index proposed by \cite{antolini_time-dependent_2005}, adapted to non-proportional hazard models, is chosen in this study. 

Consider $n$ individuals, and for $1\leq i\leq n$, $T_i$ their observation times (either survival or censoring times) and $\delta_i$ their censorship indicators. 
For $i,j = 1, \ldots,n\, i\ne j,$ we define the indicators: 
\begin{equation}
comp_{ij} = \mathbb{1}_{\{(T_i < T_j; \delta_i = 1) \cup (T_i = T_j; \delta_i = 1, \delta_j = 0)\}} \nonumber    
\end{equation}
and
\begin{equation}
    conc_{ij}^{td} = \mathbb{1}_{\{S(T_i | X_{i.}) < S(T_j| X_{j.})\}} comp_{ij}, \nonumber
\end{equation}
The estimate of the time-dependent C-index for survival models 
is equal to: 
\begin{equation}
\widehat{C}_{td} = \frac{\sum_{i=1}^n \sum_{j \ne i} conc_{ij}^{td}}{\sum_{i=1}^n \sum_{j \ne i} comp_{ij}}
\label{eq:ctd}    
\end{equation}
If we are in the proportional hazards or linear transformation models' case, the metric $\widehat{C}_{td}$ of the equation (\ref{eq:ctd}) is equivalent to the usual C-index \cite{gerdsestimating2013}.

\subsubsection{Integrated Brier Score}
The Brier score measures the squared error between the indicator function of surviving at time $t$, $\mathbb{1}_{\{T_i\ge t\}}$, and its prediction by the model $\widehat{S}(t|X_{i.})$. \cite{graf_assessment_1999} adapted the Brier score~\cite{brier_verification_1950} for censored survival data using the inverse probability of censoring weights (IPCW) and \cite{gerds_consistent_2006} subsequently proposed a consistent estimator of the Brier score in the presence of censored data.
The Brier score is defined by: 
\begin{equation}
     BS(t, \widehat{S}) = \mathbb{E}\left[ \left( Y_i (t) - \widehat{S}(t|X_{i.} )\right)^2 \right],
\end{equation}
where $Y_i(t) = \mathbb{1}_{\{T_i \ge t\}}$ is the status of individual $i$ at time $t$ and $\widehat{S}(t|X_{i.} )$ is the predicted survival probability at time $t$ for individual $i$. Unlike the C-index, a lower value of this score shows a better predictive ability of the model. 

As mentioned above, \cite{gerds_consistent_2006} gave an estimate of the Brier score in the presence of censored survival data. 
The estimate of the Brier score under right censoring is:
\begin{equation}
\widehat{BS}(t,\widehat{S}) = \frac{1}{n} \sum_{i=1}^n \widehat{W}_i(t)(Y_i(t) - \widehat{S}(t|X_{i.}))^2, \label{eq:BS}    
\end{equation}
with $n$ the number of individuals in the test set. Moreover, in the presence of censored data it is necessary to adjust the score by weighting it by the inverse probability of censoring weights (IPCW). This weighting is defined by:
\begin{equation}
    \widehat{W}_i(t) = \frac{(1-Y_i (t))\delta_i}{\widehat{G}(T_{i}|X_{i.})} + \frac{Y_i (t)}{\widehat{G}(t|X_{i.})}, \label{eq:Weights}
\end{equation}
where $\delta_i$ is the censored indicator equal to 1 if we observe the survival time and equal to 0 if the survival time is censored, and $\widehat{G}(t|x)$ is the Kaplan-Meier~\cite{kaplan_nonparametric_1958} estimator of the censored time survival function at time $t.$

The integrated Brier score~\cite{mogensen_evaluating_2012} summarizes the predictive performance estimated by the Brier score~\cite{brier_verification_1950}: 
\begin{equation}
    \widehat{IBS} = \frac{1}{\tau} \int_0^{\tau} \widehat{BS}(t, \widehat{S})dt, \label{eq:IBS}
\end{equation}
where $\widehat{BS}(t, \widehat{S})$ is the estimated Brier score and $\tau > 0$. We take $\tau > 0$ as the maximum of the observed times and the Brier score is averaged over the interval $[0,\tau].$ 
As for the Brier score, a lower value of the $IBS$ indicates a better predictive ability of the model.

\section{Results}
\label{sec:results}

In this section, we compare the performances of the Cox model with Lasso (denoted CoxL1) \cite{tibshiranilasso1997}, the neural network based on the Cox partial log-likelihood Cox-nnet \cite{chingcox-nnet:2018} presented in Section~\ref{sec:Models}, and discrete-time neural networks (NNsurv and NNsurv-deep), adapted from \cite{biganzolifeed1998} and also presented in Section~\ref{sec:Models}. The performances are compared on simulated data (with the Cox and AH models, and for several parametric configurations) and on a real data case presented below. The discrete-time C-index ($C_{td}$) and Integrated Brier Score ($IBS$) are used for this purpose. We can calculate the reference $C_{td}$ and $IBS$ values from our simulations based on the exact model used for the simulation. Note however, that the models under comparison can sometimes "beat" these reference values by chance (due to the random generation of survival times).


\subsection{Simulation study}
$n$ is the number of samples, $n\in{200;1000}$, and $p$ is the number of covariates, $p\in{10; 100; 1000}$. Note that even if our objective is to apply our models to predict survival from RNA-seq data, we present simulation results up to $1000$ covariates (instead of the potential several tens of thousands usually available with RNA-seq). Indeed, when we performed tests with 10,000 inputs, none of the model were able to perform well, thus underlining the necessity of a preliminary filtering as classically done when handling RNA-seq data \citep{conesa2016}.

\subsubsection{Results for the Cox - Weibull simulation}
The Cox-Weibull simulation corresponds to a Cox model's data with a baseline risk modeled by a Weibull distribution. In this simulation, the model satisfies the proportional hazards assumption.
The results of this simulation in \textsc{Table}~\ref{tab_CoxWeib} show that Cox-nnet performs best concerning the $C_{td}$ in all settings (regardless of the number of variables or sample size) and most settings for the $IBS$. The best $IBS$ values for Cox-nnet, as we can see from \textsc{Table}~\ref{tab_CoxWeib}, are for sample size equal to 200 and number of variables to 10 and 100 or sample size worth to 1000 and number of variables is to 100 and 1000.
CoxL1 also has the best $IBS$ (\textit{i.e.} the lowest) for a sample size of 1000 and 10 variables. These good results of CoxL1 and Cox-nnet are not surprising because we simulated the data from a Cox model. We can observe in \textsc{Table}~\ref{tab_CoxWeib} that NNsurv-deep obtains the lowest $IBS$ value for 200 individuals and 1000 variables. We can also see that the $IBS$ values of NNsurv and NNsurv-deep are very close to the reference $IBS$ values. This phenomenon is also true when the sample size is equal to 1000, and the number of variables is equal to 100.
Moreover, we can observe in \textsc{Table}~\ref{tab_CoxWeib} that some of the values of $C_{td}$ obtained for NNsurv and NNsurv-deep are close to those of Cox-nnet. We notice notably this case when the sample size is equal to 200, and the number of variables is equal to 10, and when the number of samples is 1000 and the number of variables is of 10 and 100. We can see that some of the values of the $C_{td}$ for the discrete-time neural networks are better than those obtained from the Cox model, for example, for a sample size equal to 200 and number of variables worth to 100 or for a sample size worth to 1000 and whatever the number of variables is. 

\begin{table*}[!h]
 \begin{center}
   \tabcolsep = 1.5\tabcolsep
   \begin{tabular}{p{1.5cm}|p{1.5cm}||p{1.4cm}p{1.4cm}p{1.4cm}|p{1.4cm}p{1.4cm}p{1.4cm}|}
   \hline\hline
               & n & \multicolumn{3}{|c|}{200} & \multicolumn{3}{|c|}{1000} \\
   \hline
Method & p & 10 & 100 & 1000 & 10 & 100 & 1000\\
   \hline
   \hline
Reference & $C_{td}^{\star}$ &\textbf{ 0.7442} & \textbf{0.7428}  & \textbf{0.7309} & \textbf{0.7442} & \textbf{0.7428} & \textbf{0.7309 }\\
  \cline{2-8}     & $IBS$$^{\star}$ & \textbf{ 0.0471} & \textbf{0.0549}  & \textbf{0.0582} & \textbf{0.0471} & \textbf{0.0549} & \textbf{0.0582} \\
    \hline
    \hline
NNsurv   & $C_{td}$ & 0.7137 & 0.6224  & 0.5036 & 0.7398 & 0.7282& 0.5700 \\
   \cline{2-8} &  $IBS$ &  0.0980 & 0.0646 & 0.1359 & 0.0759 & 0.0537 & 0.1007 \\
\hline
\hline
NNsurv   & $C_{td}$ & 0.7225 & 0.5982  & 0.5054 & 0.7424 & 0.7236 & 0.5741\\
  \cline{2-8}  deep &  $IBS$ &  0.0878 & 0.0689 & \textbf{0.1080} & 0.0591 & 0.0555 & 0.1185 \\
\hline
\hline
Cox   & $C_{td}$ & \textbf{0.7313}  & \textbf{0.6481} & \textbf{0.5351} & \textbf{0.7427} & \textbf{0.7309} & \textbf{0.6110} \\
   \cline{2-8}  -nnet &  $IBS$ &  \textbf{0.0688} & \textbf{0.0622} & 0.1402 & 0.0640 & \textbf{0.0498} & \textbf{0.0710} \\
   \hline
     \hline
CoxL1 & $C_{td}$ & 0.7292  & 0.5330 & 0.5011 & 0.7419 & 0.7243 & 0.5 \\
  \cline{2-8}   & $IBS$ &  0.0715 & 0.0672 & 0.1175 & \textbf{0.0541} & 0.0509 & 0.0770 \\
\hline
\hline
   \end{tabular}
\caption{Results of predicting methods on Cox-Weibull simulation} 
\label{tab_CoxWeib}
 \end{center}
\end{table*}

\paragraph{Synthesis: }
Not surprisingly, Cox-nnet has the best results on this dataset simulated from a Cox model with a Weibull distribution. However, the neural networks based on a discrete-time model (NNsurv and NNsurv-deep) have very comparable performances, and clearly outperforms the CoxL1 model when the number of variables increases.

\subsubsection{Results for the AH - Log-Normal simulation}
The results presented in \textsc{Table}~\ref{tab_simAHLN} are those obtained on the AH simulation with the baseline hazard following a log-normal distribution. In this simulation, the risks are not proportional, and the survival functions of different individuals can cross. 

We can observe that the neural networks based on a discrete-time model have the best performances concerning the $C_{td}$ and the $IBS$, and their values are close to the reference $C_{td}$ and $IBS$. This phenomenon is particularly correct for the $IBS$ when the sample size is equal to 1000, the $IBS$ values of NNsurv and NNsurv-deep are lower than those of the reference $IBS$. 
On the other hand, the methods based on the Cox partial likelihood have the highest $C_{td}$ values for a small sample size (n=200) and a small number of variables (p=10) or, on the contrary, for a large sample size (n=1000) and a large number of variables (p=1000). 
For a sample size equal to 200, neural networks based on a discrete-time model have higher $C_{td}$ values than those obtained by CoxL1 and Cox-nnet.
The values obtained for the $IBS$ by the two methods using the Cox partial likelihood are good. For a small number of individuals (n=200), the $IBS$ values of CoxL1 and Cox-nnet are very high. For example, Cox-nnet obtains $IBS$ values equal to 0.2243 and 0.1609 respectively for 10 and 100 variables, and CoxL1 gets $IBS$ values equal to 0.2278 and 0.1614, respectively. These values are very high compared to the baseline $IBS$. CoxL1 and Cox-nnet, therefore, have more difficulty with a small number of samples. The predictions of these two methods are not as good as those given by discrete-time neural networks.

\begin{table*}[!h]
 \begin{center}
   \tabcolsep = 1.5\tabcolsep
   \begin{tabular}{p{1.5cm}|p{1.5cm}||p{1.4cm}p{1.4cm}p{1.4cm}|p{1.4cm}p{1.4cm}p{1.4cm}|}
   \hline\hline
               & n & \multicolumn{3}{|c|}{200} & \multicolumn{3}{|c|}{1000} \\
   \hline
Method & p & 10 & 100 & 1000 & 10 & 100 & 1000\\
   \hline
   \hline
Reference & $C_{td}^{\star}$ & \textbf{0.7225} & \textbf{0.6857}  & \textbf{0.7070} & \textbf{0.7225} & \textbf{0.6867} & \textbf{0.7070} \\
  \cline{2-8}     & $IBS$$^{\star}$ & \textbf{ 0.0755} & \textbf{0.0316}  & \textbf{0.0651} & \textbf{0.0755} & \textbf{0.0316}  & \textbf{0.0651} \\
    \hline
    \hline
NNsurv   & $C_{td}$ & 0.6863 & \textbf{0.5971}  & \textbf{0.5358} & 0.7084 & 0.6088 & 0.5654 \\
  \cline{2-8}   &  $IBS$ &  0.1247 & \textbf{0.0780} & \textbf{0.0859} & 0.0699  & 0.0347 & 0.0533 \\
\hline
\hline
NNsurv   & $C_{td}$ & 0.7042 & 0.5793  & 0.5325 & \textbf{0.7155} &  \textbf{0.6450} & 0.5702 \\
  \cline{2-8}  deep &  $IBS$ &  0.1789 & 0.2529 & 0.1554 &  \textbf{0.0602} &  \textbf{0.0303} &  \textbf{0.0484} \\
\hline
\hline
Cox   & $C_{td}$ & \textbf{0.7128}  & 0.5812 & 0.5356 & 0.7097 & 0.6047 & \textbf{0.5720} \\
   \cline{2-8}  -nnet &  $IBS$ &  0.1342 & 0.2243 & 0.1609 & 0.0843 & 0.0875 & 0.0553 \\
   \hline
     \hline
CoxL1 & $C_{td}$ & 0.7042  & 0.5219 & 0.5112 & 0.7088 & 0.5597 & 0.5 \\
  \cline{2-8}   & $IBS$ & 0.1350 & 0.2278 & 0.1614 & 0.0608  & 0.0408 & 0.0553 \\
\hline
\hline
   \end{tabular}
\caption{Results of predicting methods on AH/Log-normal simulation} 
\label{tab_simAHLN}
 \end{center}
\end{table*}

\paragraph{Synthesis: }
On the dataset simulated from an AH model with a log-normal distribution, neural networks based on the discrete-time model have the best performances in most situations. The deep version of the model is also better than the one with only one hidden layer. In this simulation, the data do not check the proportional hazards assumption, and survival curves exhibit complex patterns for which the more versatile NNsurv-deep appears more adapted.

\subsection{Application on real datasets}

\subsubsection{Breast cancer dataset}
\paragraph{Description of data: }
The METABRIC data (for Molecular Taxonomy of Breast Cancer International Consortiulm) \cite{curtis_genomic_2012} include 2509 patients with early breast cancer. These data are available at \url{https://www.synapse.org/#!Synapse:syn1688369/wiki/27311}. Survival time, clinical variables, and expression data were present for 1981 patients, with six clinical variables (age, tumor size, hormone therapy, chemotherapy, tumor grades), and 863 genes (pre-filtered). The percentage of censored individuals is high, equal to $55\%$.

\paragraph{Results: }
The comparison results of the METABRIC dataset are summarized in \textsc{Table}~\ref{tab:metabric}. NNsurv-deep manages to get the highest value of $C_{td}$. The $C_{td}$ of NNsurv is equivalent to that of Cox, but Cox-nnet has a lower value. 
The integrated Brier score is very close for NNsurv-deep, Cox-nnet, and CoxL1, although the latter has the lowest $IBS$ value. 

On this real dataset, the differences between the models are not striking, despite the small superiority of NNsurv-deep.

\begin{table}[!h]
    \centering
    \begin{tabular}{l|c|c|c|c|c|c|}
        \hline
         & & CoxL1 &Cox-nnet & NNsurv-deep & NNsurv   \\
         \hline\hline
         Metabric & $C_{td}$ & 0.6757 & 0.6676 & \textbf{0.6853} & 0.6728 \\
        
          \cline{2-6}& $IBS$ & \textbf{0.1937} & 0.1965 & 0.1972 & 0.2038 \\
          \hline\hline
    \end{tabular}
    \caption{Results of different methods on the breast dataset (METABRIC)}
    \label{tab:metabric}
\end{table}

\section{Discussions}
This work is a study of neural networks for the prediction of survival in high-dimension. 
In this context, usual methods such as the estimation in a Cox model with the Cox partial likelihood can no longer be performed. Several methods (such as dimension reduction or machine learning methods, like Random Survival Forests \cite{ishwaran2008}) have been proposed, but our interest in this study has been directed towards neural networks and their potential for survival analysis from RNAseq data.

Two neural-network based approaches have been proposed. The first one is based on the Cox model but introduces a neural network for risk determination \citep{faraggineural1995}. The second approach is based on a discrete-time model \cite{biganzolifeed1998} and its adaptation to the high-dimensional setting was the main contribution of our work. 
In section \ref{sec:results}, we compared the standard Cox model with Lasso penalty and a neural network based on the Cox model (Cox-nnet) with those based on a discrete-time model adapted to the high dimension (NNsurv, and NNsurv-deep). To evaluate this comparison rigorously, we created a design of simulations. We simulated data from different models (Cox, AH, and AFT in appendix) with varying number of variables and sample sizes, allowing diverse levels of complexity. 

We concluded from this study that the best neural network in most situations is Cox-nnet. It can handle nonlinear effects as well as interactions. However, the neural network based on discrete-time modeling, which directly predicts the hazard risk, with several hidden layers (NNsurv-deep), has shown its superiority in the most complex situations, especially in the presence of non-proportional risks and intersecting survival curves.
On the Metabric data, NNsurv-deep performs the best, but only marginally better than the Cox partial log-likelihood-based Lasso estimation procedure, suggesting slight non-linearity and interactions.

The neural networks seem to be interesting methods to predict survival in high-dimension and, in particular, 
in the presence of complex data. The effect of censoring in these models was not studied in this work, but \cite{roblin2020use} evaluated several methods to cope with censoring in neural networks models for survival analysis. For practical applications, a disadvantage of neural networks is the interpretation difficulty. On the contrary, the output of a Cox model associated with the Lasso procedure is easily interpretable. The Cox model is therefore privileged by the domain's users nowadays. The interpretability issue of neural networks is more and more studied \citep{hao2019interpretable} and is an exciting research avenue to explore.

%
%
%
%
%
%
\bibliographystyle{ACM-Reference-Format}
\bibliography{references}  

\appendix

\section{Appendix: supplementary results}
\label{appendix}

\subsection{Simulation from the AFT - Log-Normal model}
\paragraph{\textbf{Survival times and baseline function: }}
 To simulate the data from the AFT/Log-normal model, we relied on \cite{leemis_variate_1990}. We chose to perform this simulation to generate survival data that do not respect the proportional hazards assumption. 
For this simulation, we consider that the survival times follow a log-normal distribution $\mathcal{LN}(\mu,\sigma).$ In this case, the inverse of the cumulative hazard function is expressed as (\ref{eq:InvHazLN}).
Survival times can therefore be simulated from: 
\begin{equation}
    T = \frac{1}{\exp(\beta^T X_{i.})} \exp(\sigma \phi^{-1}(U) + \mu). \label{eq:Taft} %
\end{equation}

\paragraph{\textbf{Choice of parameters of Log-Normal distribution:}}

We wish the distribution of the simulated data is close to the real data. We follow the same approach to choose the parameters $\sigma$ and $\mu$ of the survival time distribution as for the Cox/Weibull simulation presented above.  The value of the parameters is obtained from the explicit formulas: 
\begin{equation}
\mu = \ln(\mathrm{E}(T))- \frac{1}{2} \sigma^2 \label{eq:mu} \,\;\, \mbox{and} \,\;\,
\sigma^2 = \ln\left(1+\frac{\mathrm{Var}(T)}{(\mathrm{E}(T))^2}\right).
\end{equation}
Given the expectation and the standard deviation are respectively $2325$ and $1304$, the values of $\mu$ and $\sigma$ used for the simulation of the survival data should be $\mu = 7.73$ and $\sigma = 0.1760.$


\subsection{Simulation study}
\subsubsection{Results for the AFT - Log-normal simulation}
This section presents the results for data simulated from an AFT model with a baseline risk modeled by a log-normal distribution. 
The specificity of these simulated data is that they do not satisfy the proportional hazards assumption, but the survival curves do not cross.

\textsc{Table}~\ref{tab_AFTLN} shows that CoxL1 and Cox-nnet have the best results in most configurations considering $C_{td}$ or $IBS$. This good result for $C_{td}$ is particularly right when the sample size is equal to 200 or when the sample size is equal to 1000, and the number of variables is equal to 10 and 100. The $C_{td}$ obtained by the CoxL1 model is equal to 0.9867 for 200 individuals and ten variables, and the $C_{td}$ obtained for the Cox-nnet model is equal to 0.9060 for 1000 individuals and 100 variables. We can see in \textsc{Table}~\ref{tab_AFTLN} that the $C_{td}$ obtained for the neural networks based on a discrete-time model is very close to those obtained by CoxL1 and Cox-nnet and is either higher than the reference one or slightly below. For example, for a sample size equal to 200 and a number of variables equal to 10, the $C_{td}$ of NNsurv is equal to 0.9832, that of Cox-nnet is equal to 0.9867, and the reference one is equal to 0.9203. We have the same behavior for 100 variables and the same sample size or 100 variables and a sample size of 1000. 

Moreover, the $IBS$ values are the lowest for the methods based on Cox modeling in most situations. But the $IBS$ values for NNsurv and NNsurv-deep are also excellent. They are lower than the reference $IBS$ in many cases and are very close to CoxL1 and Cox-nnet. We can observe these results when the number of variables is less than or equal to 100 regardless of the sample size. 
The good results of CoxL1 and Cox-nnet might seem surprising, but we can explain it because we simulate these data from an AFT model whose survival curves do not cross. 
A method based on a Cox model will predict survival functions that do not cross. 
For this simulation, the survival function prediction obtained by CoxL1 and Cox-nnet is not cross and is undoubtedly closer to the survival function of the AFT simulation compared to discrete-time neural networks.

\begin{table}[!h]
 \begin{center}
   \tabcolsep = 1.5\tabcolsep
   \begin{tabular}{p{1.5cm}|p{1.5cm}||p{1.4cm}p{1.4cm}p{1.4cm}|p{1.4cm}p{1.4cm}p{1.4cm}|}
   \hline\hline
               & n & \multicolumn{3}{|c|}{200} & \multicolumn{3}{|c|}{1000} \\
   \hline
Method & p & 10 & 100 & 1000 & 10 & 100 & 1000\\
   \hline
   \hline
Reference & $C_{td}^{\star}$ & \textbf{0.9203} & \textbf{0.9136}  & \textbf{0.9037} & \textbf{0.9203} & \textbf{0.9136} & \textbf{0.9037} \\
  \cline{2-8}     & $IBS$$^{\star}$ & \textbf{ 0.0504} & \textbf{0.0604}  & \textbf{0.0417} & \textbf{ 0.0504} & \textbf{0.0604}  & \textbf{0.0417}\\
    \hline
   \hline
NNsurv   & $C_{td}$ & 0.9832 & 0.8349  & 0.5425 & 0.9851 & 0.9038 & 0.7426 \\
  \cline{2-8}   &  $IBS$ &  0.0265 & \textbf{0.0560} & 0.2577& 0.0247 & 0.0188 & 0.0642 \\
\hline
\hline
NNsurv   & $C_{td}$ & 0.9786 & 0.8275  & 0.5576 & 0.9857 & \textbf{0.9060}& \textbf{0.7500} \\
  \cline{2-8}  deep &  $IBS$ &  0.0295 & 0.0561 & 0.1886 & 0.0261 & 0.0207 & \textbf{0.0631} \\
\hline
\hline
Cox & $C_{td}$ & 0.9825  & \textbf{0.8558} & \textbf{0.5979} & 0.9844 & \textbf{0.9060} & 0.7085 \\
   \cline{2-8}  -nnet &  $IBS$ & \textbf{0.0122} & 0.0906 & \textbf{0.0959} & 0.0126 & 0.0374 & 0.0808 \\
   \hline
     \hline
CoxL1 & $C_{td}$ & \textbf{0.9867}  & 0.7827 & 0.5091 & 0.9856 & 0.9028 & 0.5349 \\
  \cline{2-8}   & $IBS$ & 0.0146 & 0.0965 & 0.0960 & \textbf{0.0077}  & \textbf{0.0182} & 0.0827 \\
\hline
\hline
  \end{tabular}
\caption{Results of predicting methods on AFT/Log-normal simulation} 
\label{tab_AFTLN}
 \end{center}
\end{table}
\paragraph{Synthesis: }
For data simulated from an AFT model with a log-normal distribution, Cox-nnet is the neural network with the best results in most situations when the sample size is small. When the sample size increases, NNsurv-deep is the best model considering the $C_{td}$ in most situations. Moreover, NNsurv and NNsurv-deep also seem to perform well when the number of variables is less than or equal to 100. We assume that the good results of Cox-nnet are due to the low level of complexity of the data. Indeed, the survival curves of the individuals in this dataset never cross.  

\end{document}